\documentstyle[leqno,11pt]{article}

\newtheorem{theorem}{Theorem}[section]

\def\N{{\bf N}}

\def\1{{\bf 1}}

\begin{document}

\title{On a class of arithmetic convolutions involving arbitrary sets
of integers}
\author{{\sc L\'aszl\'o T\'oth}\\
University of P\'ecs,\\
Institute of Mathematics and Informatics,\\
Hungary, 7624 P\'ecs, Ifj\'us\'ag u. 6,\\
E-mail: {\tt ltoth@ttk.pte.hu}}
\date{\sl{Math. Pannonica} {\bf 13} (2002), 249-263}
\maketitle

\begin{abstract}
Let $d,n$ be positive integers and $S$ be an arbitrary set of positive
integers. We say that $d$ is an $S$-divisor
of $n$ if $d|n$ and gcd $(d,n/d)\in S$. Consider the $S$-convolution
of arithmetical functions given by (1.1), where the sum is extended over the
$S$-divisors of $n$.

We determine the sets $S$ such that the $S$-convolution is associative and
preserves the multiplicativity of functions, respectively, and discuss other
basic properties of it.
We give asymptotic formulae with error terms for the functions
$\sigma_S(n)$ and $\tau_S(n)$,
representing the sum and the number of $S$-divisors of
$n$, respectively, for an arbitrary $S$.
We improve the remainder terms of these formulae and find
the maximal orders of $\sigma_S(n)$ and $\tau_S(n)$  assuming additional
properties of $S$.
These results generalize, unify and sharpen previous ones.

We also pose some problems concerning these topics.
\end{abstract}



\section{Introduction}

Let $\N$ denote the set of positive integers and let $S$ be an arbitrary subset
of $\N$. For $n,d\in \N$ we say that $d$ is an $S$-divisor of $n$ if $d|n$ and
gcd $(d,n/d)\in S$, notation $d|_S n$. Consider the $S$-convolution of
arithmetical functions $f$ and $g$ defined by
$$
(f *_S g)(n)=\sum_{d|_S n} f(d)g(n/d) =
\sum_{d|n} \rho_S((d,n/d)) f(d)g(n/d),  \leqno (1.1)
$$
where $\rho_S$ stands for the characteristic function of $S$.

Let $\tau_S(n)$ and $\sigma_S(n)$ denote the number and the sum of
$S$-divisors of $n$, respectively.

For $S=\N$ we obtain the Dirichlet convolution and the familiar functions
$\tau(n)$ and $\sigma(n)$. For $S=\{1\}$ we have the unitary convolution and
the functions $\tau^*(n)$ and $\sigma^*(n)$.
These have been studied extensively in the literature, see for example
\cite{McC86} and its bibliography.

Among other special cases we mention here the following ones.

Let $P$ be an arbitrary subset of the primes $p$ and $S$ be the
multiplicative semigroup generated by $P\cup \{1\}$, i. e. $S=(P)\equiv
\{1\} \cup \{n>1: p|n \Rightarrow p\in P\}$. Then the $(P)$-convolution
is the concept of the cross-convolution, see \cite{T97i}, which is a special
regular convolution of Narkiewicz-type \cite{Nar63}.

If $S$ is the set of $k$-free integers, $k\ge 2$, i. e. $S=Q_k\equiv
\{1\} \cup \{n>1: p|n \Rightarrow p^k \not| n\}$,
then the $Q_k$-divisors are the
$k$-ary divisors and (1.1) is the $k$-ary convolution, see \cite{Sur68},
\cite{Sur71}.

Let $L_k$ denote the set of $k$-full integers, i. e.
$L_k\equiv \{1\} \cup \{n>1: p|n \Rightarrow p^k| n \}$, where
$k\in \N, k\ge 2$. The $L_k$-convolution given by
$$
(f *_{L_k} g)(n)=\sum_{d|n \atop (d,n/d) \in L_k} f(d)g(n/d)  \leqno (1.2)
$$
seems to not have been investigated till now.

The aim of this note is to study some basic properties of the $S$-convolution,
to give asymptotic formulae for the functions $\sigma_S(n)$ and $\tau_S(n)$
and to investigate the maximal orders of these functions.

Assuming that $1\in S$
(then $1|_S n$ and $n|_S n$ for every $n\in \N$), we determine in Section 2 the
subsets $S$ such that the $S$-convolution is associative  and preserves
the multiplicativity of functions, respectively.

The most interesting property is that of associativity.
It turns out that, for example, the $Q_k$-convolution with $k\ge 2$ is not
associative, but the $L_k$-convolution is associative.

The $L_k$-convolution has also other nice
properties, which are analogous to those of the
Dirichlet convolution and of the unitary convolution. For example,
the set of all complex valued arithmetical functions $f$ with $f(1)\ne 0$
forms a commutative group under the $L_k$-convolution and the set of all nonzero
multiplicative functions forms a subgroup of this group.

Furthermore, let $\mu_k$ denote the inverse with respect the $L_k$-convolution
of the constant 1 function. We call it "$k$-full M\"obius function", which is
multiplicative and for every prime power $p^a$,
$\mu_k(p^a)= -1$ for $1\le a <2k$ and $\mu_k(p^a)=\mu_k(p^{a-1})-\mu_k(p^{a-k})$
for $a \ge 2k$.

Note that $\mu_1 \equiv \mu$ is the ordinary M\"obius function.
The function $\mu_2$ takes the values $-1,0,1$.

We pose the following problems: Which are the values taken by $\mu_k$?
Investigate asymptotic properties of $\mu_k$.

Note that the $S$-convolution is contained in the concept of the
$K$-convolution to be
defined in Section 2. Although there exist characterizations of basic
properties of $K$-convolutions, see \cite{Dav66} and \cite{McC86}, Chapter 4,
no study of (1.1) has been made in the literature.

Section 3 contains certain identities showing that for every $S$ the $S$-convolution
of two completely multiplicative functions can be expressed with the aid
of their Dirichlet convolution and their unitary convolution, respectively.

Asymptotic formulae with error terms for the functions $\sigma_S(n)$ and $\tau_S(n)$,
involving arbitrary subsets $S$, are given in Section 4. We show that
the remainder terms can be sharpened assuming additional properties of $S$.

In Section 5 we determine the maximal order of $\sigma_S(n)$ assuming that $S$
is multiplicative, i. e. $1\in S$ and $\rho_S$ is multitplicative, and give the
maximal order of $\tau_S(n)$ for an arbitrary $S$ with $1\in S$.

What can be said on the maximal order of $\sigma_S(n)$ for an arbitrary subset
$S$ ?

The results of Sections 4 and 5 are obtained by elementary methods,
they generalize, unify and improve the corresponding known results concerning
the functions $\sigma(n)$, $\tau(n)$, their unitary
analogues $\sigma^*(n)$, $\tau^*(n)$, those involving $k$-ary divisors and
the functions $\sigma_A(n)$, $\tau_A(n)$ associated with cross-convolutions,
see \cite{McC86}, \cite{Sur68}, \cite{Sur71}, \cite{T97i},
\cite{T97ii}.

\section{Properties of the $S$-convolution}

It is immediate that the $S$-convolution is commutative and distributive
with respect ordinary addition for every $S$.

Assume in this section that $1\in S$. Then $1|_S n$ and $n|_S n$ for every
$n\in \N$ and denoting $\delta \equiv \rho_{\{1\}}$, i. e. $\delta(1)=1$ and
$\delta(n)=0$ for $n>1$, we have $f*_S \delta=f$ for every function $f$.
This means that $\delta$ is the identity element for $*_S$.

We say that $S$ is multiplicative if $1\in S$ and
its characteristic function $\rho_S$ is multiplicative.

The $K$-convolution of arithmetical functions $f$ and $g$ is given by
$$
(f *_K g)(n)= \sum_{d|n} K(n,d) f(d)g(n/d),  \leqno (2.1)
$$
where $K$ is a complex valued function
defined on the set of all ordered pairs $\langle n, d \rangle$ with
$n,d\in \N$ and $d|n$.

For $K(n,d)=\rho_S((d,n/d))$ (2.1) becomes (1.1), therefore the $S$-convolution
is a special $K$-convolution.

\begin{theorem}
The $S$-convolution preserves the multiplicativity of functions if and only
if $S$ is multiplicative.
\end{theorem}
\vskip2mm

{\bf Proof.} It is known (\cite{McC86}, Chapter 4) that the $K$-convolution
preserves the multiplicativity if and only if
$$
K(mn,de)=K(m,d)K(n,e)
$$
holds for every $m,n,d,e\in \N$ such that $(m,n)=1$ and $d|m, e|n$.

Hence the $S$-convolution has this property if and only if
$$
\rho_S((de,mn/de))=\rho_S((d,m/d)) \rho_S((e,n/e))   \leqno (2.2)
$$
for every $m,n,d,e\in \N$ with $(m,n)=1$ and $d|m, e|n$.

If $S$ is multiplicative, then for every $m,n,d,e$ given as above
$(d,m/d)$ and $(e,n/e)$ are relatively prime,
$(de,mn/de)=(d,m/d)(e,n/e)$ and we obtain (2.2).

Conversely, if (2.2) holds and $M,N\in \N$, $(M,N)=1$ are given integers,
then taking $d=M, m=M^2, e=N, n=N^2$ we obtain
$$
\rho_S(MN)=\rho_S(M)\rho_S(N),
$$
showing that $S$ is multiplicative. $\diamondsuit$
\vskip2mm

{\bf Remark.} It follows that all the convolutions mentioned in the
Introduction preserve the multiplicativity.
\vskip2mm

\begin{theorem}
The $S$-convolution is associative if and only
if the following conditions hold:

{\rm (1)} $S$ is multiplicative,

{\rm (2)} for every prime $p$ and for every $j\in \N$ if $p^j \in S$, then
$p^{\ell} \in S$ for every $\ell >j$.
\end{theorem}
\vskip2mm

{\bf Remark.} Condition (2) is equivalent with the following:
for every prime $p$ one of the next statements is true:

(i) $p^j\in S$ for every $j\in \N$,

(ii) $p^j\notin S$ for every $j\in \N$,

(iii) there exists $e=e(p)\in \N$ depending on $p$ such that $p^j\notin S$ for
every $1\le j < e$ and $p^j\in S$ for every $j\ge e$.
\vskip2mm

{\bf Proof.} It is known (\cite{McC86}, Chapter 4) that the $K$-convolution
is associative if and only if
$$
K(n,d)K(d,e)=K(n,e)K(n/e,d/e)
$$
holds for every $n,d,e\in \N$ with $d|n, e|d$.

Therefore the $S$-convolution is associative if and only if
$$
\rho_S((d,n/d))\rho_S((e,d/e))=\rho_S((e,n/e))\rho_S((d/e,n/d))  \leqno (2.3)
$$
for every $n,d,e\in \N$ with $d|n, e|d$.

First we show that if $*_S$ is associative, then $\rho_S$ is multiplicative.
Suppose that (2.3) is satisfied, let $M,N\in \N$, $(M,N)=1$ and take
$n=M^2N^2, d=MN, e=M$. Then we have
$$
\rho_S((MN,MN))\rho_S((M,N))=\rho_S((M,MN^2))\rho_S((N,MN)),
$$
hence
$$
\rho_S(MN)=\rho_S(M)\rho_S(N).
$$
Assume now that $S$ is multiplicative. Then, taking $n=p^a, d=p^b, e=p^c$,
(2.3) is equivalent to
$$
\rho_S((p^b,p^{a-b}))\rho_S((p^c,p^{b-c}))=\rho_S((p^c,p^{a-c}))
\rho_S((p^{b-c},p^{a-b}))  \leqno (2.4)
$$
for every prime $p$ and for every $0\le c\le b \le a$. Note that it is
sufficient to require (2.4) for every $0<c<b<a$.

Suppose that $p^j\in S$, where $j\in \N$ and let $\ell >j$. We show that
$p^{\ell}\in S$.

Case 1. $\ell < 2j$. Take $a=\ell +2j, b=\ell +j, c=\ell$. From (2.4)
we obtain
$$
\rho_S((p^{\ell+j},p^j))\rho_S((p^{\ell},p^j))=\rho_S((p^{\ell},p^{2j}))
\rho_S((p^j,p^j)),
$$
$$
\rho_S(p^j)\rho_S(p^j)=\rho_S(p^{\ell})\rho_S(p^j),
$$
giving $\rho_S(p^{\ell})=1$.

Case 2. $\ell \ge 2j$. Now let $a=2\ell, b=\ell, c=\ell-j$. From (2.4)
we have
$$
\rho_S((p^{\ell},p^{\ell}))\rho_S((p^{\ell-j},p^j))=
\rho_S((p^{\ell-j},p^{\ell+j}))\rho_S((p^j,p^{\ell})),
$$
$$
\rho_S(p^{\ell})\rho_S(p^j)=\rho_S(p^{\ell-j})\rho_S(p^j),
$$
thus
$$
\rho_S(p^{\ell})=\rho_S(p^{\ell-j}). \leqno (2.5)
$$
If $\ell = kj+r$, where $k\ge 2$ and $0\le r < j$, then applying (2.5)
we have
$$
\rho_S(p^{\ell})=\rho_S(p^{\ell-j})=\rho_S(p^{\ell-2j}) = ... =
\rho_S(p^{j+r}) = 1,
$$
where $j\le j+r < 2j$ and we use the result of Case 1.

In order to complete the proof we show that if $S$ is multiplicative
and condition (2) holds, then we have (2.4) for every $0<c<b<a$.

Cosider the cases of the Remark of above. For (i) and (ii)
(2.4) holds trivially.
In case (iii) if $p^j\notin S$ for every
$1\le j \le e-1$ and $p^j\in S$ for every $j\ge e$, then (2.4) means
that the statements "[($b\ge e$ and $a-b\ge e$) and ($c\ge e$ and $b-c\ge e$)]"
and "[($c\ge e$ and $a-c\ge e$) and ($b-c\ge e$ and $a-b\ge e$)]" are
equivalent.
A quick check shows that this is true. $\diamondsuit$
\vskip2mm

{\bf Remark.} From Theorem 2.2 we obtain that the $Q_k$-convolution ($k\ge 2$)
is not associative, but the $L_k$-convolution and the (P)-convolution defined
in the Introduction are associative.
\vskip2mm

\begin{theorem}
If conditions {\rm (1)} and {\rm (2)} of {\rm Theorem 2.2} hold, then
the set of all complex valued arithmetical functions forms a commutative
(and associative) ring with identity with respect to ordinary addition and
$S$-convolution (in particular $L_k$ convolution).

This ring has no divisors of zero if and only if $S=\N$, i. e. $*_S$ is the
Dirichlet convolution.
\end{theorem}
\vskip2mm

{\bf Proof.} The first part of this result follows at once from Theorem
2.2 and from the previous remarks.

Furthermore, it is well-known that for the Dirichlet convolution there
are no divisors of zero. Conversely, suppose that $S\ne \N$
satisfies conditions (1) and (2) of Theorem 2.2. Then there exists a prime
$p$ such that $p\notin S$ and the following functions are divisors of zero:
$$
f(n)=g(n)= \cases{1, &if $n=p$,\cr
                  0, &otherwise. \cr}
\diamondsuit
$$
\vskip2mm

\begin{theorem}
If conditions {\rm (1)} and {\rm (2)} of {\rm Theorem 2.2} hold, then
the set of all complex valued arithmetical functions $f$ with $f(1)\ne 0$
forms a commutative group under $S$-convolution (in particular
$L_k$-convolution) and the set of all nonzero
multiplicative functions forms a subgroup of this group.
\end{theorem}
\vskip2mm

{\bf Proof.} This yields in a similar manner as in case
of the Dirichlet convolution and unitary convolution or in general for certain
$K$-convolutions, see \cite{McC86}, Ch.4. $\diamondsuit$
\vskip2mm

Consider now the "$k$-full"-convolution corresponding to $S=L_k$, the set
of $k$-full numbers. Let $\mu_k$ denote the "$k$-full M\"obius function",
representing the inverse of the function $I(n)=1, n\in N$ with respect to this
convolution. According to Theorem 2.4
$\mu_k$ is multiplicative and a short computation shows that
for every prime power $p^a$,
$$
\mu_k(p^a)= -1, \quad 1\le a <2k \quad
{\rm and }\quad  \mu_k(p^a) = \mu_k(p^{a-1}) - \mu_k(p^{a-k}),
\quad a \ge 2k.
$$

Observe that $\mu_1 \equiv \mu$ is the ordinary M\"obius function.

For the "squarefull M\"obius function" $\mu_2$ (case $k=2$)
we have $\mu_2(p)= \mu_2(p^2)= \mu_2(p^3)= -1$ and
$$
\mu_2(p^a) = \mu_2(p^{a-1}) - \mu_2(p^{a-2}), \quad a\ge 4.
$$

Therefore, $\mu_2(p)=\mu_2(p^2)=\mu_2(p^3)=-1,
\mu_2(p^4)=0, \mu_2(p^5)= \mu_2(p^6)=1, \mu_2(p^7)=0,
\mu_2(p^8)= \mu_2(p^9)=-1, \mu_2(p^{10})=0,...$ .

The values taken by $\mu_2$ are $-1,0,1$. This is not true for $\mu_3$,
since $\mu_3(p^a)=-1$ for $1\le a\le 5$, $\mu_3(p^6)=0, \mu_3(p^7)=1,
\mu_3(p^8)=\mu_3(p^9)=2, \mu_3(p^{10})=1, \mu_3(p^{11})=-1,
\mu_3(p^{12})=-3, \mu_3(p^{13})=-4,...$ .

We pose the following problems: Which are the values taken by $\mu_k$?
Investigate asymptotic properties of $\mu_k$.
Does it posses a mean value?

\section{Identities}

For an arbitrary $S\subseteq \N$ let $\mu_{S}$ be the M\"obius
function of $S$ defined by
$$
\sum_{d|n}\mu_{S}(n)= \rho_S(n), \quad n\in \N, \leqno (3.1)
$$
see \cite{Coh59}, therefore, by M\"obius inversion,
$$
\mu_{S}(n)=\sum_{d|n} \rho_S(d)\mu(n/d), \quad n\in \N, \leqno (3.2)
$$
where $\mu \equiv \mu_{\{1\}}$ is the ordinary M\"obius function.

The zeta function $\zeta_S$ is defined by
$$
\zeta_S(z)= \sum_{n=1}^{\infty} \frac{\rho_S(n)}{n^z}.
$$
It follows that $\zeta_{\N} \equiv \zeta$ is the Riemann zeta function
and
$$
\sum_{n=1}^{\infty} \frac{\mu_S(n)}{n^z} =\frac{\zeta_S(z)}{\zeta(z)}
\quad (z>1).  \leqno (3.3)
$$

\begin{theorem}
If $S\subseteq \N$ and $f$ and $g$ are completely multiplicative functions,
then for every $n\in \N$,
$$
(f*_S g)(n) =\sum_{d^2|n} \mu_{S}(d)f(d)g(d) (f*g)(n/d^2), \leqno (3.4)
$$
where $* \equiv *_{\N}$ is the Dirichlet convolution and
$$
(f*_S g)(n) =\sum_{d^2|n} \rho_{S}(d)f(d)g(d) (f\times g)(n/d^2), \leqno (3.5)
$$
where $\times \equiv *_{\{1\}}$ is the unitary convolution.
\end{theorem}
\vskip2mm

{\bf Proof.} Using (3.1) we have for every $n\in \N$,
$$
(f*_S g)(n) = \sum_{de=n}\rho_S((d,e))f(d)g(e) =
\sum_{de=n} \left(\sum_{j|(d,e)} \mu_{S}(j)\right) f(d)g(e).
$$
Hence with $d=ja, e=jb$,
$$
(f*_S g)(n) =\sum_{j^2ab=n} \mu_{S}(j)f(ja)g(jb)=
\sum_{j^2ab=n} \mu_{S}(j)f(j)f(a)g(j)g(b)=
$$ $$
=\sum_{j^2\ell=n} \mu_{S}(j)f(j)g(j) \sum_{ab=\ell} f(a) g(b)
=\sum_{j^2\ell=n} \mu_{S}(j)f(j)g(j) (f*g)(\ell),
$$
which is (3.4).

Furthermore,
$$
(f*_S g)(n) = \sum_{de=n}\rho_S((d,e))f(d)g(e) =
\sum_{a\in S} \sum_{de=n \atop (d,e)=a} f(d)g(e)=
$$ $$
=\sum_{a} \rho_S(a) \sum_{de=n \atop (d/a,e/a)=1} f(d)g(e).
$$
With $d=ai, e=bj$ we get
$$
(f*_S g)(n) = \sum_{a^2ij=n \atop (i,j)=1} \rho_S(a) f(a)g(a)f(i)g(j)
=\sum_{a^2b=n} \rho_S(a) f(a)g(a) \sum_{ij=b \atop (i,j)=1} f(i)g(j)=
$$ $$
=\sum_{a^2b=n} \rho_S(a) f(a)g(a) (f\times g)(b),
$$
giving (3.5).  $\diamondsuit$
\vskip2mm

\begin{theorem}
If $S\subseteq \N$, then for every $n\in \N$,
$$
\tau_S(n) = \sum_{d^2|n} \mu_{S}(d)\tau(n/d^2)
=\sum_{d^2|n} \rho_{S}(d)\tau^*(n/d^2),  \leqno (3.6)
$$
$$
\sigma_S(n) = \sum_{d^2|n} \mu_{S}(d)d\sigma(n/d^2)
=\sum_{d^2|n} \rho_{S}(d)d\sigma^*(n/d^2). \leqno (3.7)
$$
\end{theorem}
\vskip2mm

{\bf Proof.} This yields at once from Theorem 3.1 applied for $f(n)=g(n)=1$
and $f(n)=n, g(n)=1$, respectively.  $\diamondsuit$
\vskip2mm

Note that if $S$ is multiplicative, then the functions
$\tau_S(n)$ and $\sigma_S(n)$ are also multiplicative.

The generalized Euler function $\phi_S(n)= \# \{k\in \N: k\le n,
(k,n)\in S\}$ was considered in \cite{Coh59} and one has $\phi_S=\mu_S * E
=\rho_S * \phi$, where $E(n)=n, n\in \N$ and $\phi \equiv \phi_{\{1\}}$
is the ordinary Euler function, see also \cite{T97i}.

\section{Asymptotic formulae}

The following asymptotic formulae generalize and improve the known
formulae concerning the functions $\sigma(n)$, $\tau(n)$, their unitary
analogues, those involving $k$-ary divisors and the functions $\sigma_A(n)$,
$\tau_A(n)$ associated with cross-convolutions, cf. \cite{McC86}, Ch. 6;
\cite{Sur68}, Corollary 3.1.1; \cite{Sur71}, Corollary 3.1; \cite{T97i},
Theorem 12; \cite{T97ii}, Theorem 2; see also \cite{T98}, Corollary 1.

\begin{theorem}
If $S\subseteq \N$, then
$$
\sum_{n\le x} \sigma_S(n) =\frac{\zeta(2) \zeta_S(3)}{2\zeta(3)}x^2 + R_S(x),
\leqno (4.1)
$$
where the remainder term can be evaluated as follows:

{\rm (1)} $R_S(x)=O(x\log^{8/3} x)$ for an arbitrary $S$,

{\rm (2)} $R_S(x)=O(x\log^{5/3} x)$ for an $S$ such that
$\sum_{n\in S} \frac1{n}< \infty$ (in particular for every finite $S$)
and for every multiplicative $S$,

{\rm (3)} $R_S(x)=O(x\log^{2/3} x)$ for every multiplicative $S$ such that
$\sum_{p\notin S} \frac1{p}< \infty$ (in particular if the set
$\{p: p\notin S\}$ is finite).
\end{theorem}
\vskip2mm

{\bf Proof.} We have from (3.7),
$$
\sum_{n\le x} \sigma_S(n) = \sum_{d\le \sqrt{x}} \mu_S(d)d \sum_{e\le x/d^2}
\sigma(e).
$$
Applying now the well-known result of Walfisz \cite{W63},
$$
\sum_{n\le x} \sigma(n) =\frac{\zeta(2)}{2}x^2 + O(x\log^{2/3} x)
$$
we obtain
$$
\sum_{n\le x} \sigma_S(n) = \sum_{d\le \sqrt{x}} \mu_S(d)d
\left(\frac{\zeta(2)x^2}{2d^4} + O\left(\frac{x}{d^2} (\log
\frac{x}{d^2})^{2/3}\right) \right)=
$$ $$
= \frac{\zeta(2)x^2}{2} \sum_{d=1}^{\infty} \frac{\mu_S(d)}{d^3}+
O\left(x^2 \sum_{d>\sqrt{x}} \frac{|\mu_S(d)|}{d^3}\right)+
O\left(x (\log x)^{2/3} \sum_{d\le \sqrt{x}} \frac{|\mu_S(d)|}{d}\right).
$$
For the main term apply (3.3) and the given error term yields from the next
statements:

(a) For an arbitrary $S\subseteq \N$, $|\mu_S(n)|\le \sum_{d|n} \rho_S(d)
\le \tau(n)$ for every $n\in \N$ and
$$
\sum_{n\le x} \frac{|\mu_S(n)|}{n}\le \sum_{d\le x} \frac{\rho_S(d)}{d}
\sum_{e\le x/d} \frac1{e}= $$ $$= O\left(\log x \sum_{d\le x}
\frac{\rho_S(d)}{d}\right)
= \cases{O(\log x), &if $\sum_{n=1}^{\infty} \frac{\rho_S(n)}{n} < \infty$,\cr
         O(\log^2 x), &otherwise. \cr}
$$

(b) If $S$ is multiplicative, then $\mu_S$ is multiplicative too,
$\mu_S(p^a)=\rho_S(p^a)-\rho_S(p^{a-1})$ for every prime power $p^a$ ($a\ge 1$)
and $\mu_S(n)\in \{-1,0,1\}$ for each $n\in \N$.

(c) Suppose $S$ is multiplicative. Then
$$
\sum_p \sum_{k=1}^{\infty} \frac{|\mu_S(p^k)|}{p^k} \le
\sum_p \left(\frac{|\rho_S(p)-1|}{p} + \sum_{k=2}^{\infty} \frac1{p^k}\right)
= \sum_{p\in S} \frac1{p(p-1)} + \sum_{p\notin S} \frac1{p-1}\le
$$
$$
\le 2\left(\sum_{p\in S} \frac1{p^2} + \sum_{p\notin S} \frac1{p}\right) <
\infty \quad
{\rm if} \quad \sum_{p\notin S} \frac1{p} < \infty.
$$
It follows that in this case the series $\sum_{n=1}^{\infty}
\frac{|\mu_S(n)|}{n}$ is convergent.  $\diamondsuit$
\vskip2mm

\begin{theorem}
If $S$ is an arbitary subset of $\N$, then
$$
\sum_{n\le x} \tau_S(n) =\frac{\zeta_S(2)}{\zeta(2)}x \left(\log x+ 2\gamma-1
+\frac{2\zeta'_S(2)}{\zeta_S(2)} -\frac{2\zeta'(2)}{\zeta(2)} \right)
+O(\sqrt{x}\log^2 x), \leqno (4.2)
$$
where $\gamma$ is the Euler constant and $\zeta'_S(z)$ is the derivative of
$\zeta_S(z)$.
\end{theorem}
\vskip2mm

This result follows applying the first identity of (3.6) and using
Dirichlet's formula
$$
\sum_{n\le x} \tau(n) =x (\log x+ 2\gamma -1)+ O(x^{\alpha}).
$$

The remainder term of (4.2) can be improved assuming further properties
of $S$. For example, if $S$ is multiplicative, then the error term is
$O(\sqrt{x} \log x)$ and if $S$ (i. e. $\rho_S$) is completely multiplicative
and $\{p: p\not\in S\}$ is a finite set, then the error term is
$O(x^{\alpha})$. We do not go into details.

\section{Maximal orders}

Generalizing the result of Gronwall concerning the function $\sigma(n)$
we prove the following theorem.

\begin{theorem}
Let $S$ be an arbitrary multiplicative subset. Denote by $P$ the
set of primes $p$ such that $p^j\in S$ for every $j\in \N$.
For every $p\notin P$ let $s(p)\in \N$ denote the least exponent $j$
such that $p^j\notin S$ (i. e. $p^j\in S$ for every $1\le j <s(p)$ and
$p^{s(p)}\notin S$).

Then
$$
\limsup_{n\to \infty} \frac{\sigma_S(n)}{n \log \log n}=e^{\gamma}
\prod_{p\notin P} \left(1- \frac1{p^{2s(p)}}\right).
$$
\end{theorem}
\vskip2mm

{\bf Proof.} For every $p\in P, a\in
\N$ and for every $p\notin P, a < 2s(p)$
the $S$-divisors of $p^a$ are all divisors $1,p,p^2,...,p^a$. Hence
$\sigma_S(p^a)= \sigma(p^a)= 1+p+p^2+...+p^a$.

For every $p\notin P$ and $a\ge 2s(p)$ the numbers $p^{s(p)}$ and $p^{a-s(p)}$
are certainly not $S$-divisors of
$p^a$, since $(p^{a-s(p)},p^{s(p)})=p^{s(p)}
\notin S$. Therefore
$\sigma_S(p^a)< (1+p + p^2+ ... +p^{a -s(p)-1}) +
 (p^{a -s(p)+1}+...+p^a)
< p^{a -s(p)}+ p^{a -s(p)+1}+...+p^a
\le p^{a -2s(p)+1}+ p^{a -2s(p)+2}+...+p^a$.

We obtain that
$$
\frac{\sigma_S(p^a)}{p^a} \le 1 + \frac1{p} + \frac1{p^2}+...+
\frac1{p^{2s(p)-1}} \leqno (4.3)
$$
holds for every prime power $p^a$ with $p\notin P$
with equality for $a=2s(p)-1$.

Also, for every $p\in P, a\in \N$,
$$
\frac{\sigma_S(p^a)}{p^a} < \left(1-\frac1{p}\right)^{-1}.
\leqno (4.4)
$$

We show that
$$
\frac{\sigma_S(n)}{n} \le  e^{\gamma} \prod_{p\notin P}
\left(1- \frac1{p^{2s(p)}}\right)
\log \log n (1+o(1)) \qquad {\rm as }\qquad n\to \infty.
$$
Using (4.3) and (4.4) we have for every $n\ge 1$,
$$
\frac{\sigma_S(n)}{n} \le \prod_{p|n \atop{p\in P}} \left(1 -\frac1{p}
\right)^{-1}
\prod_{p|n \atop{p\notin P}} \left(1+\frac1{p} +\frac1{p^2}+...+
\frac1{p^{2s(p)-1}}\right)=
$$
$$
= \prod_{p|n \atop{p\le \log n \atop{p\in P}}} \left(1 -\frac1{p}\right)^{-1}
\prod_{p|n \atop{p> \log n \atop{p\in P}}} \left(1 -\frac1{p}\right)^{-1}
\prod_{p|n \atop{p\le \log n \atop{p\notin P}}} \left(1+\frac1{p} +\frac1{p^2}+...+
\frac1{p^{2s(p)-1}}\right) \times $$ $$ \times
\prod_{p|n \atop{p> \log n \atop{p\notin P}}} \left(1+\frac1{p} +\frac1{p^2}+...+
\frac1{p^{2s(p)-1}}\right) \le
$$
$$
\le \prod_{p\le \log n \atop{p\in P}} \left(1 -\frac1{p}\right)^{-1}
\prod_{p\le \log n \atop{p\notin P}} \left(1+\frac1{p} +\frac1{p^2}+...+
\frac1{p^{2s(p)-1}}\right) \times $$ $$ \times
\prod_{p|n \atop{p> \log n \atop{p\in P}}} \left(1 -\frac1{p}\right)^{-1}
\prod_{p|n \atop{p> \log n \atop{p\notin P}}} \left(1 -\frac1{p}\right)^{-1}=
$$
$$
=\prod_{p\le \log n \atop{p\notin P}} \left(1 -\frac1{p^{2s(p)}}\right)
\prod_{p\le \log n} \left(1 -\frac1{p}\right)^{-1}
\prod_{p|n \atop{p>\log n}} \left(1 -\frac1{p}\right)^{-1}\le
$$
$$
\le \prod_{p\le \log n \atop{p\notin P}} \left(1 -\frac1{p^{2s(p)}}\right)
\prod_{p\le \log n} \left(1 -\frac1{p}\right)^{-1}
\prod_{p|n \atop{p>\log n}} \left(1 -\frac1{\log n}\right)^{-1}=
$$
$$
=e^{\gamma} \prod_{p\notin P} \left(1- \frac1{p^{2s(p)}}\right)
\log \log n (1+o(1)),
$$
applying Mertens' theorem
$\prod_{p\le x} (1-\frac1{p})=\frac{e^{-\gamma}}{\log x}(1+o(1))$
as $x\to \infty$, and the fact that
$\#\{p: p|n, p>\log n\}\le \frac{\log n}{\log \log n}$.

Now we show that this upper bound is asymptotically attained.

For a given $\varepsilon >0$ choose $t$ so large such that
$$
\prod_{p>t} \left(1-\frac1{p^2} \right) \ge 1 - \varepsilon.
$$

For this $t$ choose an exponent $a\ge 1$ such that
$$
\prod_{p\le t} \left(1-\frac1{p^a} \right) \ge 1 - \varepsilon.
$$

Consider the sequence $(n_k)_{k\ge 1}$ given by
$$
n_k=\prod_{p\le t \atop{p\in P}} p^{a-1}
\prod_{p\le t \atop{p\notin P}} p^{2s(p)-1}
\prod_{t< p\le e^k} p.
$$

We obtain
$$
\frac{\sigma_S(n_k)}{n_k}=
\prod_{p\le t \atop{p\in P}} \left(1+\frac1{p} +\frac1{p^2}+...+
\frac1{p^{a-1}}\right)  \times $$ $$ \times
\prod_{p\le t \atop{p\notin P}} \left(1+\frac1{p} +\frac1{p^2}+...+
\frac1{p^{2s(p)-1}}\right)
\prod_{t< p\le e^k } \left(1+\frac1{p}\right) \ge
$$
$$
\ge \prod_{p\le t} \left(1- \frac1{p^a}\right)
\prod_{p\notin P} \left(1- \frac1{p^{2s(p)}}\right)
\prod_{p>t} \left(1-\frac1{p^2}\right)
\prod_{p\le e^k} \left(1 -\frac1{p}\right)^{-1} \ge
$$
$$
\ge (1-\varepsilon)^2 \prod_{p\notin P} \left(1- \frac1{p^{2s(p)}}\right)
e^{\gamma} k (1+o(1)) \qquad {\rm as } \qquad k\to \infty,
$$
applying Mertens' theorem again.

Furthermore, considering the Chebysev function $\theta(x)=\sum_{p\le x} \log p$
and using the elementary estimate $\theta(x) =O(x)$, we get
$$
\log n_k \le O(1) + \theta(e^k) =O(e^k).
$$
Hence, for sufficiently large $k$,
$$
\log \log n_k \le O(1) + k < (1+\varepsilon)k.
$$
Therefore
$$
\limsup_{k\to \infty} \frac{\sigma_S(n_k)}{n_k \log \log n_k} \ge
\frac{(1-\varepsilon)^2}{1+\varepsilon} e^{\gamma}
\prod_{p\notin P} \left(1- \frac1{p^{2s(p)}}\right),
$$
and the proof is complete.  $\diamondsuit$
\vskip2mm

A direct consequence of Theorem 5.1 is the following result.

\begin{theorem}
Let $S$ be an arbitarary multiplicative subset
and suppose that there exists $s\in \N$ such that for every prime $p$,
$p^j\in S$ for every $1\le j <s$ and $p^s\notin S$.
Then
$$
\limsup_{n\to \infty} \frac{\sigma_S(n)}{n \log \log n}=
\frac{e^{\gamma}}{\zeta(2s)}.
$$
\end{theorem}
\vskip2mm

This result can be applied for $S=Q_k$ (case $s=k\ge 1$), for $S=L_k$
(case $s=1$).
\vskip1mm
What is the maximal order of $\sigma_S(n)$ for an arbitrary subset $S$ ?

\begin{theorem}
Let $S$ be an arbitrary subset such that $1\in S$. Then
$$
\limsup_{n\to \infty} \frac{\log \tau_S(n)\log \log n}{\log n}=\log 2.
\leqno (4.5)
$$
\end{theorem}
\vskip2mm

{\bf Proof.} It is well-known that this result holds for the function
$\tau(n)$ (case $S=\N$) and that for the sequence $n_k=p_1p_2...p_k$,
where $p_i$ is the $i$-th prime,
$$
\lim_{k\to \infty} \frac{\log \tau(n_k)\log \log n_k}{\log n_k}=\log 2.
$$

Taking into account that if $1\in S$, then $\tau_S(n)=\tau(n)$ for every
squarefree $n$ and $\tau_S(n)\le \tau(n)$ for every $n\in \N$, (4.5)
follows at once. $\diamondsuit$

\end{document}